# Developing Optimization Models with Cognitive Systems Engineering

Tyler C. O'Brien, MS[1], Emily L. Tucker, PhD[1*], Steven Foster, PhD[1], Sudeep Hegde, PhD[1]

[1]Clemson University, Department of Industrial Engineering, Clemson, SC, USA

* Corresponding author: Emily Tucker etucke3@clemson.edu

## Abstract

One goal of applied operations research is to improve decisions in practice. This requires modelers and stakeholders to have a shared understanding of the system and for the developed model to reflect the system's core dynamics. There are four areas to address: the underlying problem must be understood, the mathematical formulation of the problem must be representative of the system at hand, the data must be appropriate, and the model-generated recommendations must be understandable by the stakeholders. While developing models, operations researchers may primarily rely on past experience in model development, rather than underlying theory, to guide decisions on how to include stakeholders in the modeling process. In parallel, the field of Cognitive Systems Engineering has developed methodologies and practices to understand systems, stakeholder needs, and environments. To improve the rigor of the "application" in applied operations research, we present a framework to integrate Cognitive Systems Engineering methods with optimization model development. We apply the integrated framework to a case study of locating hand sanitizer stations in response to COVID-19 at a large academic institution.

## 1. Introduction

Applied operations research (OR) modelers aim to develop models and methods that will support real-world decisions. The toolkit includes mathematical optimization, which allows modelers to generate recommendations quickly and for problems far larger than people can solve to optimality unassisted. Optimization has been used successfully to improve decisions in many domains, including fulfilling retail orders (Andrews et al., 2019), generating routes for trucking (Rönnqvist et al., 2017), scheduling crews for air travel (Barnhart et al., 2003), prepositioning inventory (Duran et al., 2011), generating new districts for organ transplants (Gentry et al., 2015), and scheduling sports leagues (Durán et al., 2021). These applications of optimization have led substantial improvements in cost, time, worker fulfillment, and human health.

There are necessary conditions for model-generated recommendations to be effective in practice. Modelers must have (i) a sufficient understanding of the problem, (ii) the ability to express the problem mathematically and solve it, (iii) appropriate data to parameterize the model, and (iv) the ability to communicate recommendations to decision-makers. We suggest that unless someone is embedded in an application domain, most OR modelers are only experts in (ii) and sometimes (iv). To support the remainder, modelers must work with stakeholders and decision-makers in the application domain itself.

Consider Figure 1. There are three distinct entities in the system: modelers, the model, and stakeholders. Figure 1A shows the area in which modelers are trained and tend to thrive: model development. We refer to it as the nominal view. In applied OR, though, the reality is broader and looks more like Figure 1B, i.e., the contextual view. The modelers work with stakeholders to understand the problem and the context. The model recommendations are interpreted by the stakeholders, and each entity is affected by surrounding context (Eisenberg et al., 2019). Here, we use “context” to refer to any factor that affects the people involved or decisions that could be made; this could involve incentives, institutional requirements, history, culture, and available technology. Without a proper perspective on the context and stakeholder needs, perceptions, and constraints, a modeler may develop a substandard model or make the situation worse.

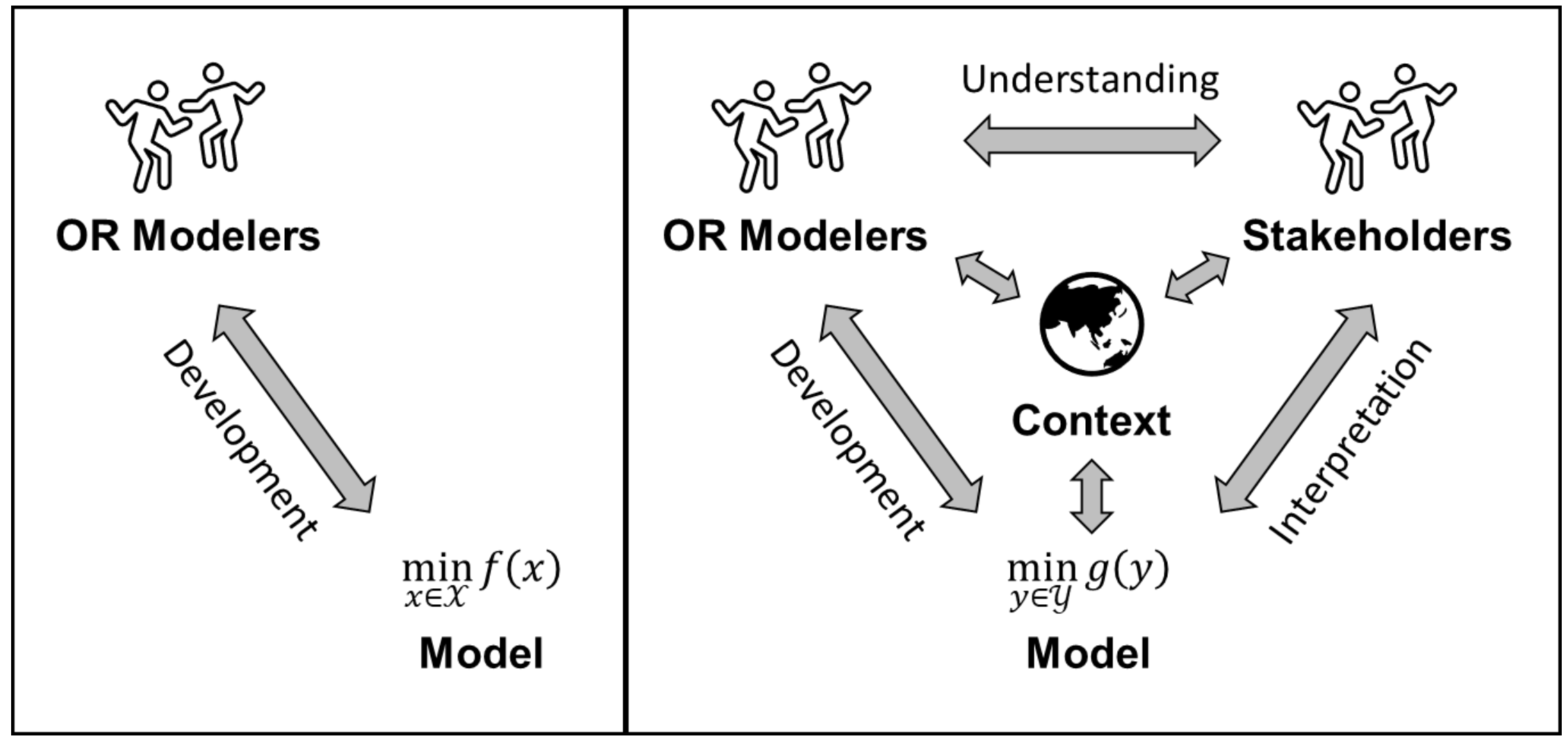

**Figure 1.** Nominal view (panel A) vs. contextual view (panel B) of models.

To move from the nominal to the contextual view, modelers often rely on their own experiences or best practices for modeling (e.g., Michael Pidd, 1999). However, there are opportunities to ground the process in existing theoretical frameworks.

In this vein, we turn to the field of Cognitive Systems Engineering (CSE) (Hollnagel et al., 2006; Hollnagel & Woods, 2005). CSE focuses on understanding the interactions between humans, the environment, technologies and processes in a system, and the properties that emerge thereof. It does so through methods to elicit perspectives of stakeholders as they engage with and adapt to a dynamic world, even as they cope with uncertainty, complexity, and various pressures such as time. CSE methods have led to numerous applications, including the design of decision interfaces, accident modeling, systems dynamics modeling, and organizational policy forming (e.g., Burns & Hajdukiewicz, 2017; Goncalves Filho et al., 2019; McGeorge et al., 2015).

Traditionally, CSE methods have focused on *qualitatively* modeling systems and the decision processes involved. These approaches translate stakeholder perspectives, interactions, and mental models of the world into insights for improving the system as a whole. A mental model is an agent's (a human or other dynamic system) internal conceptual representation of an external system, in other words, a 'view' or understanding of how a system works. Each individual's mental model is shaped by their knowledge and their interactions with the system, and influences decision-making (Doyle & Ford, 1998). Mental models can vary based on agent's dynamic interactions with the system. CSE methods can help explicate the mental models and tacit cognitive processes of domain stakeholders and develop insights to dynamically inform the design and optimization of systems. This approach, however, has rarely been integrated with mathematical modeling techniques in general and has not yet been applied in particular to the design of optimization models.

In this paper, we propose a framework to integrate CSE methods into optimization model development. Through subsequent iterations of model versions, stakeholders are involved to increase modelers' understanding of the problem and the decision-makers' understanding of the model simultaneously. The goal of our method is to create greater alignment in the mental models of the real-world-facing decision-makers with the modelers who seek to represent parts of the world through mathematical expressions. This enables the stakeholders to provide further insights to the limitations of the model and increase their understanding of the problem through the lens of the current model version.

The approach brings together CSE methods of interviewing and thematic analysis with optimization model development to support iterative model design. It is flexible in its ability to engage multiple stakeholders and encompass different modeling approaches. To illustrate its use, we present a case study of a COVID-19 mitigation strategy: where to locate hand sanitizer stations on a university campus. We use semi-structured interviews (Hegde et al., 2020) to understand initial and follow-up decisions of decision-makers and the priorities of key stakeholders. Models are iteratively developed and revised in accordance with the interview findings. The results are recommendations for a campus-wide allocation of stations.

The primary contribution of this paper is an *integrated OR-CSE framework to improve decisions.* To the best of our knowledge, this is the first study to integrate these two complementary fields. It has particular utility in settings with emerging information, where mental models of both modelers and stakeholders may evolve over time.

This paper is part of a larger project that is focused on campus adaptation due to COVID-19 and organizational resilience. One of the model versions included in this paper was published as a brief work-in-progress paper in conference proceedings (O'Brien et al., 2021). In a paper focused on the qualitative CSE analyses, we observe that the creation of a centralized structure for the stakeholder groups on campus increased adaptive capacity (Foster et al., 2023).

The remainder of the paper has the following structure. In Section 2, we discuss related literature. In Section 3, we present the OR-CSE integrated framework. We apply the framework to a case study of campus hand sanitizer station locations in Section 4. In Section 5, we discuss the framework and case study, and in Section 6, we conclude.

## 2. Literature Review

Work to improve the practice of operations research has spanned from guidance on best practices to empirical research studies. Among the perspective pieces, Michael Pidd (1999) presents six key principles of modeling which include a focus on simplicity, with model development beginning with small and growing in complexity over the process of modeling. Morris (1967) describes the "enriching" process as model development proceeds as additional detail added if the model remains tractable. Geoffrion (1976) draws attention to *insight* as the underlying goal of OR. He suggests that simple, auxiliary models that transparently capture core dynamics can be effectively used alongside more complex models. This is particularly important with emerging risks, e.g., Morgan et al. (2022), where transparency, trust, and relationships are critical to develop models as contextual information evolves.

Kunz et al. (2017) focus on modeling for humanitarian logistics and similarly emphasize the importance of knowing context for research to be able to have practical insight. Kovacs and Moshtari (2019) present a meta-process for humanitarian operations research, based on Sodhi and Tang (2014), that is initiated with exploratory research and continues in collaboration with stakeholders.

Among empirical and qualitative research, Willemain (1994) presents a survey of twelve modeling experts that finds that context, assessment, and structure are important to the modeling process. In a follow-up, Willemain (1995) studies how experts formulate models via an experiment with twelve modelers who solved exercises. Notably, the experts did not linearly progress through the stages of modeling, i.e., from context to implementation, but rather revisited stages, including context, over time. Pidd (2010) considers the question of how models are used. His interviews suggest that investigating and improving are common goals for model use. They may be developed for new or existing decision processes. Sharkey et al. (2024) categorizes modeling efforts in response to fundamental surprise events; these are events for which it is impossible to plan ahead of time, e.g., COVID-19. Based on interviews of modelers and stakeholders, they develop five categories of model adaptation, from changing the data to reflect the new context to developing a new model. Other research has studied how models may affect decision-making (Chung et al., 2000).

Murphy (2005) suggests that a "theory of practice" for operations research would necessarily include cognitive science. He presents an initial list of key skills for a practitioner, ranging from technical to business aspects. Two of the areas, cognitive problem mapping and the ability to abstract are important for modeling and understudied. The taxonomy presented provides a matrix of model types based on the model goal and different levels of clarity on the state of the problem. In this paper, we provide an approach to integrate a cognitive-focused theory, CSE, into the model development process.

Throughout, researchers underscore that working with stakeholders is crucial (Sharkey et al., 2021). de Gooyert et al. (2017) review OR papers that include stakeholders and the connections to stakeholder theory. They describe four types of problems in which stakeholders are involved, i.e., optimization models, trade-offs (termed "balancing"), problem structuring, and broadening perspectives (termed "involving"). These are relevant in any application domain, though are particularly popular in community-based OR (Johnson et al., 2018), where qualitative data may be used to incorporate stakeholder perspectives. Quantitative approaches include probability encoding (Spetzler & Stael von Holstein, 1975), where model inputs are informed by domain expertise. User interfaces can be used to communicate model-generated recommendations and/or allow stakeholders to interact with the model (Ahani & Trapp, 2021). Further, when the goal of the modeling process is to improve policy, Gass (1983) notes additional stakeholder involvement, i.e., external evaluation and public debate, is necessary, beyond standard steps of validation and assessment. Other fields have developed frameworks for participatory research that include the seminal ladder of citizen participation (Arnstein, 1969) and community-based participatory research (Jones & Wells, 2007).

Soft OR methods are often used in "messy" situations, including when there are multiple stakeholders and unclear problem definitions (Mingers, 2011; Mingers & Rosenhead, 2004). These methods include Soft Systems Methodology, Strategic Choice Approach, and Strategic Operations Development and Analysis.

Within the context of "hard" operations research, where this paper is situated, research has included specific aspects of model development. Merrick and Weyant (2019) connect concepts from information theory to model development and provide an approach to determine model granularity. Pollack and Steimle (2022) study how to aggregate states in Markov Decision Processes. Shehadeh and Tucker (2022) evaluate models with and without distributional ambiguity in disaster relief settings; they that find mis-specifying the underlying distribution may increase shortages.

Willemain et al. (2003) study model-based decision-support systems and highlight that one cause of poor outcomes is the use of the wrong model. One approach to mitigate structural uncertainty, i.e., the case when the right model to use is not known precisely, comes from the field of health economics. The modeler may conduct scenarios analyses using alternate model structures and present each to the decision-makers (Bojke et al., 2009; Strong et al., 2012). Note this is different from a Pareto-analysis approach (e.g., Daskin and Tucker, 2018) where the trade-offs between decision alternatives are generated from the same underlying (ostensibly correct) model. After development, model validation is critical; best practices are to conduct face, internal, external, and predictive validation (Eddy et al., 2012).

CSE methods are well-suited to elicit perspectives from domain experts based on their lived experience of interacting with real-world systems, and to translate findings into tangible insights that inform design and development of systems (Elm et al., 2008; Militello et al., 2010; Woods & Roth, 1988). Traditionally, these approaches have been applied in complex sociotechnical domains, such as healthcare and transportation management, to inform the design and development of interfaces that support users in monitoring complex dynamic processes and decision making (e.g., Elix & Naikar, 2021; Hettinger et al., 2017). Such methods also have the potential to be used to develop and update mathematical models based on user-driven insights and help bring greater alignment between the mathematical models and the user's mental models. However, this intersection of CSE and mathematical modeling for real-world decision-making has previously not been explored.

The case study example provided in this paper is focused on a university's response to COVID-19. Many OR researchers successfully developed models to support decisions for higher education (see brief review in Sharkey et al., 2024). These included the redesign of a campus bus network (Chen et al., 2022) and class scheduling (Barnhart et al., 2022; Navabi-Shirazi et al., 2022) and cohorting (Gore et al., 2022). OR researchers led a discrete-choice experiment to identify student preferences (Steimle et al., 2022) and collaborated closely with stakeholders (Barnhart et al., 2022). Two problem-definition reports also helped to lower the barrier to entry to develop useful optimization (Shen, 2020) and simulation (Currie et al., 2020) models.

Taken together, this paper presents a new approach to ground optimization model development within CSE theory. This deepens the literature in OR related to stakeholder-involved model development and amplifies the potential impact of CSE.

## 3. OR-CSE Integrated Modeling Approach

In this section, we present the integrated CSE-OR framework for optimization model development (Figure 2). Modelers and stakeholders iteratively update their conceptual understanding (i.e., mental models), and modelers develop models of the decision process. Each iteration represents the development of one model version. After the decision-makers consider a model version sufficiently representative, the process concludes.

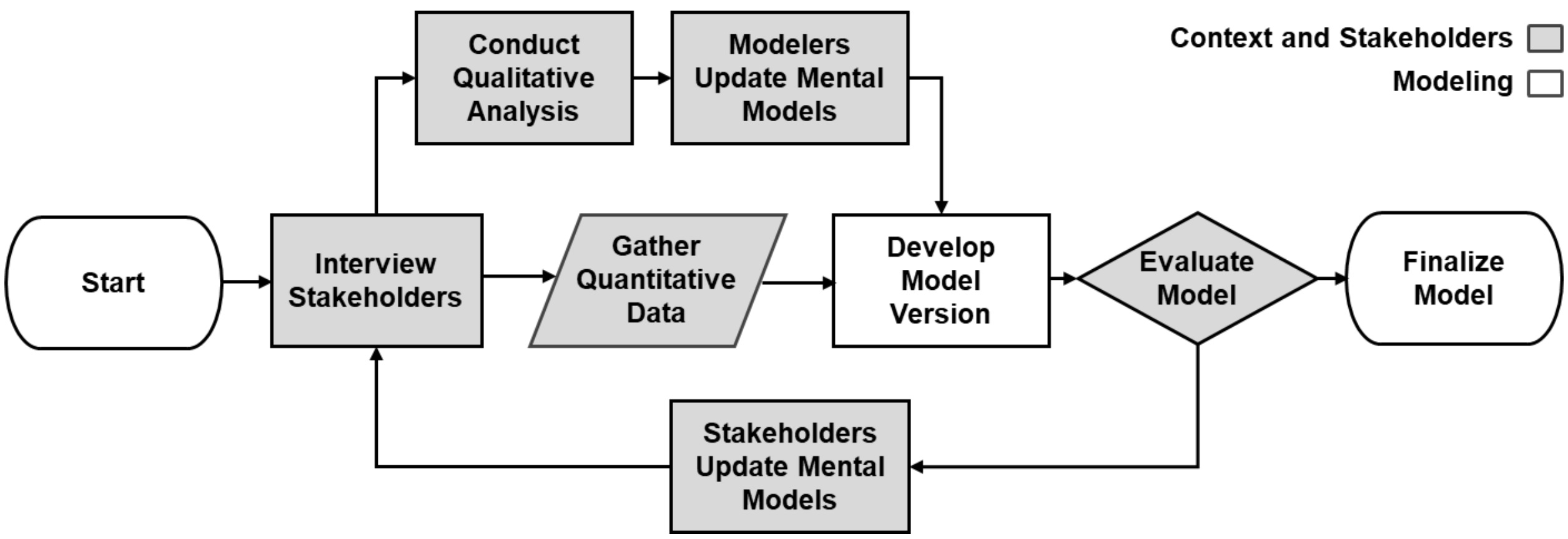

**Figure 2**. OR-CSE modeling process.

In the diagram of the process (Figure 2), the gray boxes represent steps that may change conceptual understanding, and the white boxes indicate model development (Figure 1B). We note that not every step may be necessary each iteration. The details of each step are as follows.

*Start.* The process begins when modelers and stakeholders discuss the potential for an optimization model to make recommendations for practice. This step can be initiated either by the modelers or the stakeholders.

*Interview Stakeholders.* The team interviews key stakeholders to gain insight into stakeholder perspectives and underlying context. Ideally, interviews are conducted with representatives from each stakeholder group, including the group that will make the final decisions, i.e., decision-makers, and the groups that would be affected by the model recommendations. If the decisions have the potential for widespread impact, particular care should be taken to interview stakeholders from marginalized groups and those who have less power to identify and avoid unintended consequences.

One approach to interviewing is to use a semi-structured process. Stakeholders are asked both pre-planned questions and questions that arise during the interview, including clarifications. Beforehand, the modeling team develops a semi-structured guide using CSE methods, as below.

These techniques elicit knowledge of domain experts on how decisions are made adaptively under varying situations characterized by time pressure, uncertainty, and other factors (Hegde et al., 2020; Klein et al., 1989). The guide is designed to reveal the underlying mental models of the stakeholders to better understand their perspectives and how decisions are made.

The structure of the semi-structured interview guide is a series of open-ended questions and "cognitive probes." The open-ended questions allow the participant to provide an account of the situation(s) encountered. The cognitive probes target specific aspects of the actual workflows, e.g., the types of pressures and constraints, adaptive strategies and decisions, and factors that enable or impede execution of the decisions. The interview questions are tailored to the specific context and the roles of the participants. In subsequent iterations (Figure 2), the team conducts follow-up interviews; these include further questions related to their decisions, challenges and strategies, types of data used, coordination with other stakeholder groups within and outside the organization.

*Conduct Qualitative Analysis.* The team analyzes the qualitative data (e.g., transcripts) from the interviews to identify themes and insights that are related to model development. A theme is defined in this study as a pattern of ideas or concepts identified in several excerpts across participant interviews. First, the analysts conduct qualitative coding, where segments of text are classified as related to an idea, i.e., a "code" (Mayring, 2004; Miles & Huberman, 1994). Two analysts independently code segments and then reconcile any code disagreements. Next, the thematic analysis is conducted (Braun & Clarke, 2012). Text excerpts from interview are grouped based on shared ideas (codes), and analysts derive themes that occur across stakeholders.

*Modelers Update Mental Models.* The modelers consider the results from the qualitative analyses, including the thematic analyses and generated insights, and revise their understanding of the system, accordingly. This understanding may or may not be directly applicable to the next model iteration. The primary goal at this step is to improve the alignment of the modelers' and stakeholders' perceptions of the external system and decision process.

*Gather Quantitative Data.* In parallel, relevant quantitative data are collected. Additional sources for data may have been identified through the interviews.

*Develop Model Version.* In this step, a new or revised version of the model is created. It reflects the modelers' current understanding of the context and stakeholder perspectives. New or revised model-relevant features or feedback that were identified in the qualitative analysis are incorporated. These may include revisions of decisions, constraints, objectives, parameters, and/or uncertainty (cf. categorization in Sharkey et al. 2024).

Large revisions to the design may be necessary to reflect a new understanding of priorities and restrictions, e.g., the modelers may switch the framework from a network flow optimization to integer program. As iterations progress, the changes may become smaller, e.g., revising the objective function or constraints. New solution methods may also be developed at this step, as appropriate.

*Evaluate Model Iteration.* The current model version is assessed by the modelers and stakeholders to determine if it is sufficient to make recommendations for practice. Assessment may take different forms. The process could include analyses with current data, sensitivity analyses, inventory of model assumptions, and/or formal model validation (Eddy et al., 2012). To support the evaluation, the modelers may share visuals and/or data with the decision-makers. These materials could include model recommendations (e.g., optimal solution), performance metrics (e.g., optimal value, results by subgroup), and assumptions. The model recommendations may be compared with historical decisions, if available. If the model is determined by decision-makers to be sufficiently representative of practice, the process concludes. If not, the process continues.

*Stakeholders Update Mental Models.* During the evaluation of the current model version, stakeholders may realize that there are aspects of the system of which they had not been aware. The differences between their previous mental models and their new awareness of the system may come directly from model insights or as an indirect response to seeing it abstracted in a new way. These updates could include the recognition of an issue of which they were unaware, a new strategy, or the importance of a particular dynamic. They then implicitly update their mental models of the external system. Stakeholders may also update their mental models as they learn more about the context, e.g., during in emerging crises, or as the context changes over time.

*Finalize Model.* If decision-makers are satisfied with the model during the evaluation step, development concludes. Model analyses may be conducted to answer key policy questions and/or the model may be put in production. Over time, tactical updates may be made to the model. These may include re-running analyses using updated data or minor modeling changes. If the decision-making context changes or substantially new information is revealed, the model may be revised by re-starting the process.

## 4. Case Study: Access to Hand Sanitizer

In this section, we describe a case study that applies the integrated CSE-OR framework to develop models to distribute hand sanitizer stations across a university campus. We discuss the context of the problem in Section 4.1 and the stakeholders, interview protocol, and qualitative analyses in Section 4.2. We present the initial iterations focused on problem scoping in Section 4.3. The main iterations that develop five model versions are in Section 4.4. Key model results are presented in Section 4.5.

### 4.1. Context

The onset of the COVID-19 pandemic forced universities to shift to remote operations rapidly in the spring of 2020 (Birmingham et al., 2023). To support a safe return to campus in Fall 2020, Clemson University implemented a portfolio of mitigation strategies to reduce risks to students, staff, and faculty. One strategy was to deploy hand sanitizer stations across campus (Centers for Disease Control and Prevention, 2020). These movable stations are approximately 4 feet tall and were placed inside buildings near exterior doors (Appendix Figure A-1).

In this example, the Facilities department was the decision-maker. Stakeholders included the departments associated with the campus Emergency Operations Center (EOC) as well as other individuals across campus. The EOC was the body responsible for coordinating COVID-19 response.

### 4.2. CSE Approach: Knowledge Elicitation

Throughout the OR-CSE integrated process, there were 11 stakeholder participants (Figure 3). They represent Facilities, Procurement, Environmental Health and Safety, Student Health Services, Housing and Dining, the Provost's Office, and the EOC. These were the groups that were most closely involved in decisions related to planning, procuring, allocating, deploying, and monitoring hand sanitizer stations around campus. Their roles spanned multiple organizational levels, from leadership to boots-on-the-ground frontline work. These departments and individuals were connected both through the EOC and through pre-existing professional relationships between individual departments (Figure 3). The participants were identified using snowball sampling, as each participant was connected through the EOC. Each interview was conducted on Zoom by two or more team members and recorded for later qualitative analysis.

The study was reviewed by the Clemson University Institutional Review Board and approved under the 'exempt' category. The interview protocol, including the question guide, is provided in Appendix Section A-1.

The interviews focused on five major categories of questions. The first set of questions was related to the participant's role, specifically at the beginning of the pandemic and how it related to the university's sanitizer deployment. They were also asked to describe the difficulties that they and their related department/team faced at the start of the pandemic in the Spring of 2020. For instance, "*What were the main challenges faced at the start of the pandemic (around March 2020)?*". The next set of questions focused on information availability related to sanitizer station placement considerations and included probes such as about the information they gathered from external sources (e.g., CDC) and internal university sources related to hand sanitation policies. Third, the participants were asked about their station placement strategies. These included questions about goals and metrics and the use of data and models, such as, *"What were your main priorities with regards to hand sanitation on campus?"*. The fourth set focused on monitoring hand sanitizer use and situational changes (e.g., class modality changed from online

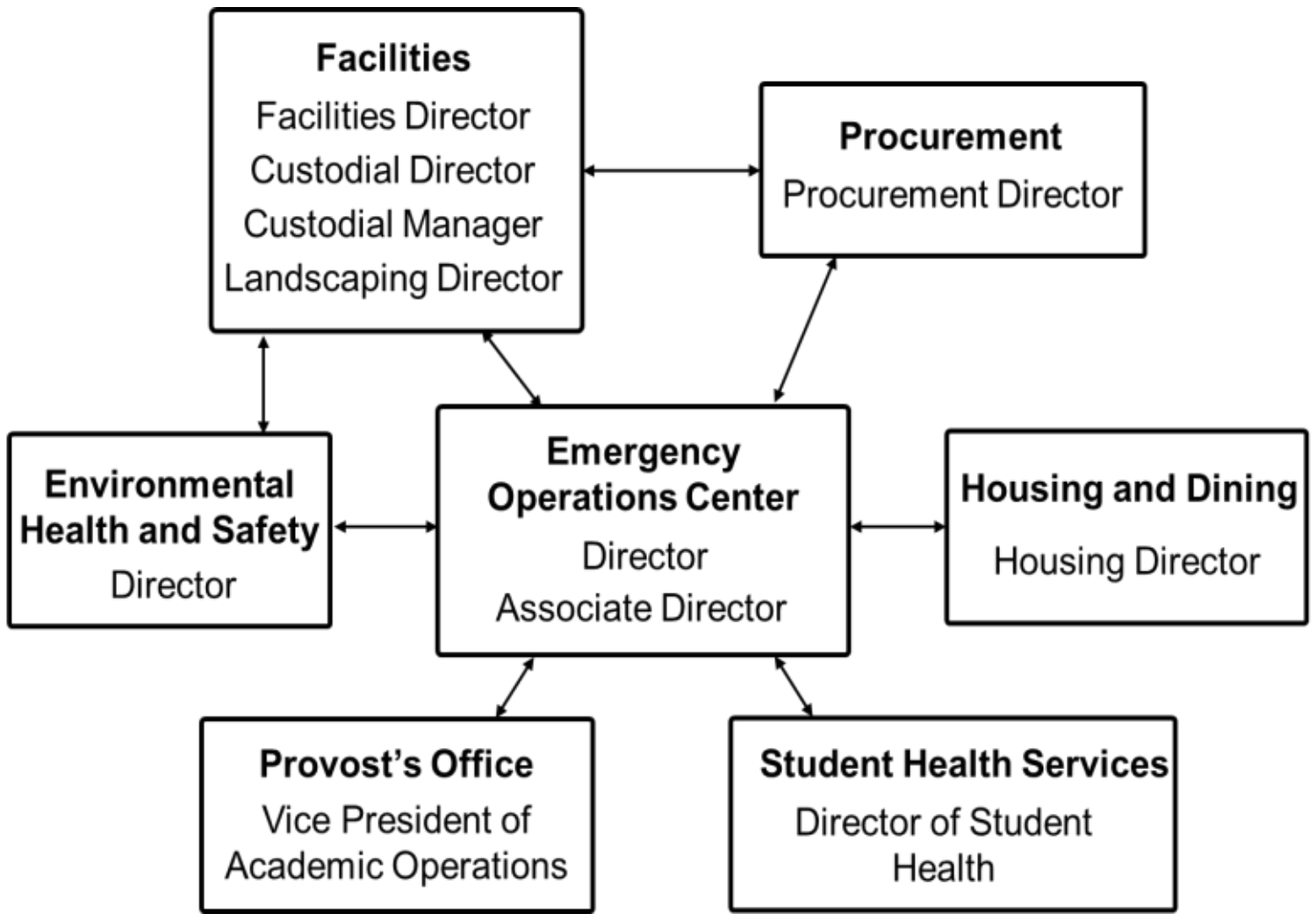

**Figure 3.** Stakeholder groups and titles of interview participants.

to a hybrid). The last set asked about directions for the future.

The interviews were semi-structured and included several probes to allow further questions related to topics of particular interest. In additional to questions that focused on improving the academic team's mental models of decision-process more generally, several questions were included to elicit information to support model design specifically. These included questions on goals, priorities, and strategies which have relevance to an optimization model's objective functions and decision variables. Following interviews, qualitative analysis via coding and subsequent thematic analysis were performed to identify key decision-making considerations of the participant stakeholders shown in Figure 3.

### 4.3. Problem Scoping

At the beginning of the project, our academic team started with a broad goal to support the Facilities department in COVID-19 response efforts. We had an interest in hand sanitizer decisions, but we did not have a specific problem definition *a priori*. The first set of interviews focused on understanding the priorities, challenges, and key decisions during the pandemic. Based on these, we developed initial ideas for decisions and model structures (without math) and conducted rapid iterations with the Facilities department. Each was based on our initial understanding of the context and system, which deepened over time. We did not include other stakeholders at this stage.

First, we considered models to recommend station locations based on the movement of students, e.g., network flow-type models such flow-capture and path-based approaches. However, the Facilities department noted that each station needed to be inside a building (due to early pandemic risk of theft and ease of replenishment), and hence, it was not necessary to consider exterior, between-building flow. It was also unlikely that we would be able to collect sufficient data to support models of within-building flow. We next considered whether to optimize station replenishment, e.g., with a routing-based approach on sanitizer station refills. They noted that this would not be as helpful, because stations were regularly checked.

The Facilities department described the process of deciding where to put stations at the beginning of the pandemic as a “guessing game” based on perceptions of foot traffic. We proposed that we could develop a model to allocate stations to buildings, and they agreed that that would be helpful. They stated that the focus should be on academic buildings as there was a different process for residence and dining halls.

During the scoping phase, the Facilities department recommended multiple data sources as potential proxies for foot traffic, i.e., demand. One was “door access control data” that was collected via scans of identification cards each time a controlled door was opened, e.g., all exterior entrances. The department that managed the identification cards provided time-stamped data on scan times and door openings for each controlled door. We anonymized the data and aggregated traffic by week; a heat map of scanned entrances per building is provided in Appendix Figure A-2. A limitation is that multiple people may enter per scanned access; this data may have underestimated true demand. An alternative data source, Wifi access points, was also suggested, but it did not become available during the project period.

### 4.4. Model Development

After the initial rapid iterations to scope the problem, we conducted five full iterations to develop models to optimize hand sanitizer deployment. In this subsection, we present the details of the progression of the OR-CSE process (Figure 2) and the model associated with each iteration. A summary of the models is in Table 1, and notation is in Table 2.

**Table 1.** Model Summary.

| | V1 | V2 | V3 | V4 | V5 |
|---|---|---|---|---|---|
| **Name** | Facilities Heuristic | Single Building Max Coverage | Target Coverage | Amplified Coverage by Door | Max Coverage by Door |
| **Scope** | Campus-wide | One building | Campus-wide | Campus-wide | Campus-wide |
| **Granularity** | Each building | Multiple interior locations | Each building | Each exterior door of each building | Each exterior door of each building |
| **Decision Variables** | Number of stations per building | Select candidate locations; classroom coverage | Number of stations per building | Number of stations per door | Select candidate door |
| **Objective** | n/a | Maximize number of covered classrooms | Minimize squared difference between demand and supply | Maximize amplified coverage of demand | Maximize coverage of demand |
| **Constraints** | Proportional to demand; at least one per building | p stations, coverage distance | At least one station per building; only allocate available stations | Pre-assigned number of stations per building; at most two stations per door | At least one station per building; at most one station per door; allocate available stations |

**Table 1**. Model Notation.

| Sets | Description | Applicable to Model Version(s) |
|---|---|---|
| $B$ | Set of buildings | 1, 3, 4, 5 |
| $C$ | Set of classrooms | 2 |
| $I$ | Set of door indices | 4, 5 |
| $J$ | Set of candidate locations within one building | 2 |
| **Parameters** | | |
| $a_{cj}$ | Coverage indicator (1 if $t_{cj} < T$, 0 otherwise for classroom $c \in C$ and candidate location $j \in J$) | 2 |
| $d_b$ | Demand in building $b \in B$ | 1, 3 |
| $\hat{d}_{bi}$ | Demand in building $b \in B$ at door $i \in I$ | 4, 5 |
| $l_{bi}$ | Door indicator (1 if building $b \in B$ has door index $i \in I$, 0 otherwise) | 4, 5 |
| $N$ | Total number of stations available | 1, 3, 5 |
| $p$ | Number of stations to be located | 2 |
| $s_b$ | Number of stations pre-allocated to building $b \in B$ | 4 |
| $T$ | Coverage threshold time | 2 |
| $t_{cj}$ | Time to walk from classroom $c \in C$ to candidate location $j \in J$ | 2 |
| $u$ | Number of available uses (pumps) per station | 3 |
| **Decision Variables** | | |
| $x_b$ | Number of stations located in building $b \in B$ (integer) | 1, 3 |
| $\hat{x}_{bi}$ | Number of stations located in building $b \in B$ at door $i \in I$ (integer) | 4 |
| $\bar{x}_{bi}$ | 1 if station is located in building $b \in B$ at door $i \in I$, 0 otherwise (binary) | 5 |
| $y_j$ | 1 if station is located at candidate $j \in J$, 0 otherwise (binary) | 2 |
| $z_c$ | 1 if classroom $c \in C$ is covered, 0 otherwise (binary) | 2 |

#### 4.4.1. Model Version 1 (Heuristic)

At the outset, the academic team's mental models were as follows: we understood that (i) the primary goal was to distribute stations across campus, (ii) the focus was on academic buildings rather than residence halls and dining, and (iii) each building should receive at least one station. We had access to door access control data, and exterior door openings were considered a reasonable proxy for foot traffic entering each door.

Based on this understanding, we developed Model V1, a heuristic baseline to represent the process of the initial, ad hoc allocation of stations to buildings.

$$x_b = \left\lceil \frac{d_b}{\sum_{b' \in B} d_{b'}} * N \right\rceil \quad \forall b \in B \tag{1}$$

The heuristic, given in (1), decides the number of stations $x_b$ to place in each building $b \in B$ out of $N$ stations available based on the proportion of the building's demand ($d_b$) vs. demand across campus. Demand is defined as the number of times an exterior door in building $b \in B$ is accessed. Higher traffic buildings receive proportionally more stations, and each building receives at least one.

During the evaluation of Model V1, we made the following observations. It does not guarantee an optimal solution to the general problem of distribution and does not specify locations within buildings. By rounding the number of assigned stations up to the nearest integer, the heuristic may require a greater number of stations than available. It also assumes that allocating to locations with higher demand would result in the stations being used more frequently.

We presented the model-recommended allocation to the Facilities department (Section 4.5, Figure 4). One of their takeaways was that their actual, ad hoc distribution based on perception (i.e., not with door access control data) was reasonably close to a heuristic that used quantitative data. This validated their existing mental models of perceived demand across campus and indicated that model-supported recommendations could amplify their expert judgment. We continued model development to address its limitations.

#### 4.4.2. Model Version 2 (Single Building Max Coverage)

The Facilities department expressed interested in assigning locations within each building, and a new data source was identified and collected, i.e., schematic diagrams of each building on campus. We next developed a model to optimize station locations within a single building (Model V2).

The modelers' mental models indicated that within an academic building, classrooms would be the primary demand points. As data on user behavior, e.g., hand sanitizer use based on location, was not available, the focus should be on improving access to the stations, i.e., coverage. Based on this conceptual understanding, Model V2 is a max coverage facility location model that determines the optimal locations for stations within a single building (2).

$$\max \sum_{c\in C} z_c \qquad (2.1)$$

s.t.

$$\sum_{j\in J} a_{cj} y_j \geq z_c \quad \forall c \in C \qquad (2.2)$$

$$\sum_{j\in J} y_j = p \qquad (2.3)$$

$$y_j, z_c \in \{0,1\} \quad \forall j \in J, c \in C \qquad (2.4)$$

The objective function (2.1) maximizes the number of classrooms that are covered by a sanitizer station. Candidate locations, $J$, include doorways, hallway intersections, and other locations along interior hallways. Constraints (2.2) define whether a location is covered by at least one station, i.e., if the travel time between them is less than a threshold time designated via a binary coverage parameter. Constraint (2.3) locates exactly $p$ stations. Constraints (2.4) ensure the decision variables are binary.

In version evaluation, we noted that the model was limited by data availability. Identifying candidate station locations and the times between candidate locations and classrooms would be time-consuming and require the involvement of building staff members. By pre-fixing the number of stations per building, $p$, the model does also not allow analyses of campus-wide allocation.

Furthermore, the focus of the real-world decision-process shifted. The Facilities department acquired smaller dispensers for hand sanitizer and placed one in each classroom on campus. The framing of Model V2, i.e., locating stations to cover classrooms, was no longer necessary.

### 4.4.3. Model Version 3 (Target Coverage)

Stakeholders changed their focus to be the allocation of stations to locations near building entrances, rather than other locations in building interiors. We learned during interviews that entrances do not serve the same purpose in each campus building. For example, in Freeman Hall two entrances are primarily used, whereas in the library, there is a single, heavy-use entranceway. This would affect station placement; the library entrance may require more than one station per door.

With revised mental models, the modelers developed Model V3. This version of the model allocates stations to buildings across campus, and the assignment is based on a target coverage metric. It seeks to match the number of stations to the building with the estimated demand for sanitizer at that building. The formulation is as follows:

$$\min \sum_{b \in B} (d_b - u x_b)^2 \qquad (3.1)$$

s.t.

$$x_b \geq 1, \qquad \forall b \in B \qquad (3.2)$$

$$\sum_{b \in B} x_b \leq N \qquad (3.3)$$

$$x_b \geq 0, \text{ integer} \qquad \forall b \in B \qquad (3.4)$$

The objective function (3.1) minimizes the difference from target coverage using a quadratic function. The amount of available sanitizer is based on the number of uses (pumps) per station multiplied by the number of stations allocated to the building. Constraints (3.2) ensures that each building receives at least one hand sanitizer station. Constraint (3.3) limits the total number of located stations across all buildings to the total number of stations available. Constraints (3.4) requires the number of stations assigned to be non-negative and integer.

During evaluation, we noted that the model does not provide specific locations for the stations to be placed; this is non-trivial for buildings with many doors. Model V3 gives priority to high volume buildings rather than priority to high volume doors, which they noted was more important. The model also penalizes allocations that provide more stations than projected need; stakeholders noted that this was not a problem in practice if stations are available.

### 4.4.4. Model Version 4 (Amplified Coverage by Door)

The stakeholders indicated that one of the limitations (specific location recommendations within buildings), was important to address. We developed Model V4 to recommend placement within buildings. It aims to maximize the coverage of pedestrian traffic across campus. Eligible locations are the exterior doors, i.e., the entry and exits to each building. Doors may have up to two stations assigned to each to consider doors that are higher use. To quickly evaluate the impact of the model switch, we fix the number of stations per building $b \in B$ to a user-defined parameter, $s_b$. Demand at door $i \in I$ in building $b \in B$ is $\hat{d}_{bi}$, and integer decision variables $\hat{x}_{bi}$ determine how many stations to place at each door $i$ in building $b$. The formulation is given as:

$$\max \sum_{b \in B} \sum_{i \in I} l_{bi} \hat{d}_{bi} \hat{x}_{bi} \qquad (4.1)$$

s.t.

$$\sum_{i \in I} \hat{x}_{bi} = s_b \qquad \forall b \in B \qquad (4.2)$$

$$\hat{x}_{bi} \leq 2 \qquad \forall b \in B, i \in I \qquad (4.3)$$

$$\hat{x}_{bi} \geq 0, \text{ integer} \qquad \forall b \in B, i \in I \qquad (4.4)$$

The objective function (4.1) maximizes the amplified coverage of demand across campus. Amplified coverage is defined as the product of demand per door and the number of stations assigned per door. A parameter coefficient $l_{bi} \in \{0,1\}$ indicates whether building $b \in B$ has a door with index $i \in I$, as not all buildings have the same number of doors. Constraints (4.2) ensure that the total number of stations allocated to doors in each building is equal to a user-defined parameter, e.g., current number of stations in the building. Constraints (4.3) allow up to two stations per door. Constraints (4.4) enforce integer decisions.

One limitation is that each building is restricted to its current capacity, which means that the model does not address cross-campus allocation. Because the model maximizes traffic coverage and some doors can receive more than one station, the model automatically allocates the maximum number of stations to the highest demand doors. This is reasonable in some cases; however, if doors have similar demand, the lower demand door may not be allocated a station. The Facilities department appreciated the shift to by-door allocation and were interested in evaluating campus-wide allocation.

#### 4.4.5. Model Version 5 (Max Coverage by Door)

Model V5 identifies where to locate stations in buildings and the number of stations to place in each building. It is very similar to V4, with the following exceptions: only a single station can be allocated to each door ($\bar{x}_{bi} \in \{0,1\}, \forall b \in B, i \in I$), and the number of stations per building is not fixed. Its formulation is:

$$\max \sum_{b \in B} \sum_{i \in I} l_{bi} \hat{d}_{bi} \bar{x}_{bi} \qquad (5.1)$$

s.t.

$$\sum_{i \in I} \bar{x}_{bi} \geq 1 \qquad \forall b \in B \qquad (5.2)$$

$$\sum_{b \in B} \sum_{i \in I} \bar{x}_{bi} = N \qquad (5.3)$$

$$\bar{x}_{bi} \in \{0,1\} \qquad \forall b \in B, i \in I \qquad (5.4)$$

The objective function (5.1) maximizes the coverage of demand across campus. Constraints (5.2) ensure that each building has at least one station. Constraint (5.3) ensures that the number located does not exceed the number available (*N*). Constraints (5.4) enforces that the decision variables are binary.

We noted in model evaluation that the model accounts for exterior doors only, and it may be worthwhile to consider some high-traffic internal building locations as well, e.g., student lounges, based on additional interviews. At this stage, the Facilities department was satisfied with Model V5 and did not think further iterations would improve overall decisions related to station allocation. Based on this, Model V5 was designated the final version.

### 4.5. Model Results

We ran analyses for each model based on the 36 academic buildings on campus. These buildings represented 55% of campus foot traffic. For each model, demand was measured as foot traffic by door, or in aggregate by building, based on one week of door access data from February 2021. We also performed sensitivity analyses using three weeks of data from April 2021. The allocations with data from February and April are very similar. We present results using February 2021 data. The models were solved in Excel (V1) and in Python (V2-V5) using the Pyomo package and Gurobi as the solver.

The baseline, real-world, comparator is the original deployment by the Facilities department in Spring 2020. They placed 102 sanitizer stations across 36 academic buildings on campus. Each building received at least one station, and the average number of stations per building was 2.8. The College of Business Building received the most with 20.

Three of the model versions allocate stations to buildings across campus (V1, V3, and V5). The full allocations are shown in Appendix Table A-1. Across the models, the average of number of stations per building remains fairly stable (3.3, 2.8, and 2.8 stations, respectively vs. 2.8 in initial deployment), but the number of buildings with a single station varies considerably (6, 24, and 10 buildings, respectively vs. 15 in initial deployment). Model V3 is a target coverage model, and the goal is to match stations to foot traffic. We observe this dynamic with the low traffic buildings; even if there are many entrances, the model attempts to avoid overage penalties and assigns the minimum number of stations per building. The last version of the model (V5) is more indicative of stakeholder goals and focuses on traffic area (door) rather than aggregated by building. The intuition here is that the buildings do not experience an aggregated demand; people in a building may not all be in the same area. Rather, doors are the bottleneck areas where demand is realized.

The recommendations from the model compared with the initial deployment are presented in Figure 4. The x-axis shows the difference in the number of stations assigned to a building from the model vs. the initial deployment. For example, the bin of "[1-3]" indicates buildings where the model recommends one, two, or three more stations than were initially deployed. Model V1 recommends that 19 should have higher allocations, and nine buildings match initial deployment. Model V3 recommends that six buildings increase (e.g., Brackett Hall initially had two but is recommended to have four), and 16 buildings remain the same. Model V5 recommends 14 buildings increase, and 14 buildings remain the same.

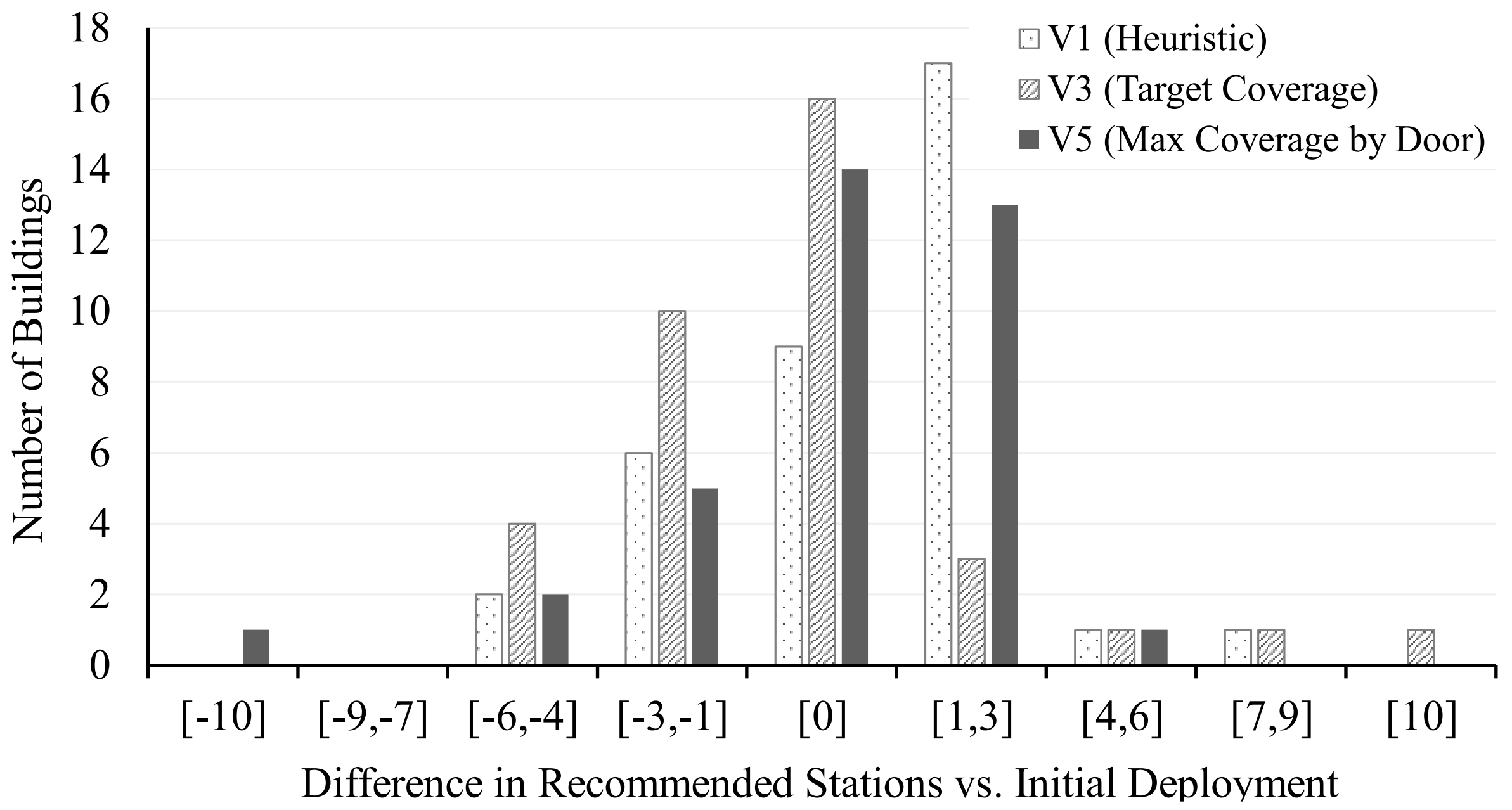

**Figure 4.** Comparison of Model-Generated Allocations to Initial Deployment.

Variations of Figure 4 were shown to the Facilities department throughout the model development process. We found that it was helpful to present the results in comparison to the initial deployment (vs. presenting them without this context). It provided a grounding of the results with a real-world setting with which they were familiar.

To demonstrate the use of the single building max coverage model (Model V2), we used Freeman Hall, where the Industrial Engineering department is located (Figure 5, Facilities Support Services, 2020). There were 10 candidate locations for stations including exterior doors, hallway intersections, and key interior locations. With four stations, the selected locations are the largest classroom ($j = 10$) and three interior intersections ($j = 1,5,9$). Recommended locations

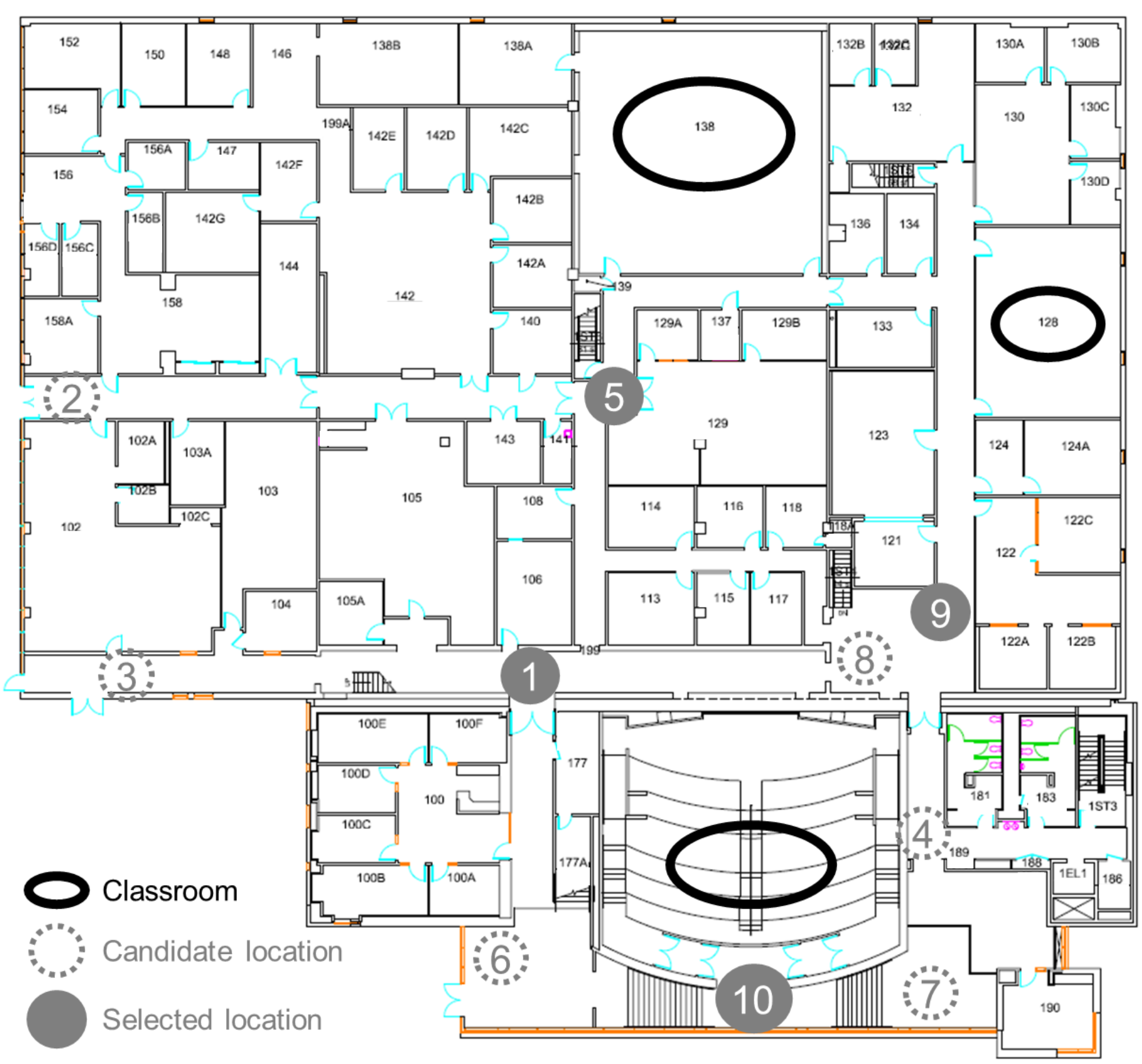

**Figure 5.** Recommended locations for four stations in Freeman Hall (Model V2).

with one to three stations allocated are shown in Appendix Table A-2.

Model V4 locates stations at exterior doors across campus, allows up to two stations to be placed at each, and limits each building to its initial allocation. For Freeman Hall, the initial allocation was three with three doors with one station each; Model V4 recommends that one door receive two stations, and another door receives one station.

## 5. Discussion

We present an integrated approach to incorporate contextual information from stakeholders into the development process of optimization models (Figure 2). The main OR contribution of this framework is to ground its foundation on the qualitative methods from the field of Cognitive Systems Engineering (CSE). CSE methods, such as knowledge elicitation of domain experts, have traditionally been used to develop insights for design of decision aids, such as visual interfaces, or work systems and processes (e.g., Bennett et al., 2023; Bisantz & Burns, 2009; Burns & Hajdukiewicz, 2017). Such studies typically focus on helping designers better understand the needs and perspectives of domain stakeholders, the ultimate users of the designed system. However, this is the first study to combine the context-rich CSE methods with the technical capabilities of OR methods to improve decision-making. The study uniquely demonstrates a method to bridge gaps in the mental models of multiple groups, including OR modelers, CSE/design researchers, and domain stakeholders. Developing a shared understanding of the system is key step toward integrating perspectives for engineering design to be effective.

The framework uses an iterative method such that the modelers develop versions based on their increasing understanding of context. It is adaptable to evolving situations, including fundamental surprise events of COVID-19. The iterative approach has the byproduct of deepening stakeholders' understanding of both optimization models and their own context as they interpret the model results. We apply the framework to a case study of hand sanitizer locations at a large public institution in the US in partnership with the Facilities department.

The process aims to increase the resulting model's usefulness to stakeholders. For a model to be useful, the modeler must understand the domain context and be able to communicate model results. The framework enables this increasing modeler understanding through on-going interactions with stakeholders. These perspectives affect the mental models that modelers develop about the system and hence the optimization models themselves.

In the presented case study, the decision-maker, the Facilities department, was the primary evaluator of each model iteration. However, in other settings, this level of interaction with the decision-maker may not be feasible, but their involvement in some capacity is critical (Ahani & Trapp, 2021), e.g., in at least one iteration. Decision-makers and stakeholders are considered partners in the process of design, and this framework connects with the theme of participatory research more broadly (Arnstein, 1969; Jones & Wells, 2007; Kovacs & Moshtari, 2019). Throughout the iterative process, the framing of stakeholders as experts on the problem is important and should not be lost. Simple phrases, such as describing model recommendations as "building on the great work already done" by decision-makers, are true and can go a long way in continuing to build trust.

There is flexibility in the framework. While the case study uses semi-structured interviews as the qualitative methodology, alternatives such as focus groups and surveys may also be useful. One caveat with surveys with close-ended questions is that they are more constrained by the perspective of the person creating the survey. If a modeler designs a close-ended survey based on an incorrect mental model of the situation, it may not lead to a better understanding of the true situation; rather, it may reinforce incorrect mental models. Different stakeholders may be interviewed at different iterations. In the hand sanitizer case study, the decision-makers (i.e., the Facilities department) were interviewed each iteration, and other of the other stakeholder groups (e.g., Dining) were interviewed only once. We note that we did not interview students as stakeholders; this is a limitation of the developed models. Not every step needs to be completed in each iteration, e.g., quantitative data collection may happen only as needed.

In the model evaluation step, the modelers present the results and/or model design to stakeholders, and they help evaluate the quality of the model version. The evaluation step may be conducted with the decision-makers (as in our case example) or with multiple stakeholder groups. They may indicate that the current model is or is not sufficient in its representation or usefulness. Transparency is important. The overall model structure, including decisions and assumptions, should be clear at each iteration. Outputs can be communicated via visuals, and providing a comparator based on their own context may be helpful. For example, this could include presenting recommendations in comparison to their experience or showing diagrams with landmarks they are familiar with. We suggest that very limited math, if at all, should be shown.

We note that the models developed in the case study are not complex but focus on the key features of the decisions. Our goal was not to eliminate all assumptions and limitations. In a crisis, a faster and good decision is more important than a slow but “optimal” one, particularly when the situation such as COVID-19 evolves quickly (Sharkey et al., 2024). The idea is to improve optimization modeling *within* the uncertainty of the real-world; it is not always necessary to capture the uncertainty within the model itself.

The study demonstrates how the convergence of CSE and OR modeling could be an effective strategy for organizations to adapt scalable models to fast-changing situations. For CSE experts, it is important to not only capture evolving real-world patterns but also to keep pace with modeling needs as they continue to elicit knowledge from the domain stakeholders. The study represents a unique challenge for CSE, to not just understand stakeholder needs but to engage with mathematical modelers and to help facilitate a three-way mental-model reconciliation, i.e. between the domain stakeholders, modelers, and CSE experts as designers.

There are open questions in how to integrate conflicting information from stakeholders and disagreements in priorities. The hand sanitizer case study is one of strong alignment; each group was focused on providing a clear public health benefit without competition for resources. Differences only arose in the perspectives of constraints and priorities. In other settings, however, decision-makers and stakeholder groups may disagree on the goal itself, e.g., how to zone land in a city. It would be worthwhile to study how to use a CSE-framework to develop multi-objective models that could support decisions in contested settings. Further research could also consider methods for long-term model maintenance and use. The framework ends with a completed model (Figure 2), but the underlying context will shift over time. A more advanced process could incorporate longer-term evaluation and additional iterations of the model. Studies could consider whether ad hoc triggers or regular reviews would be effective to restart model development iterations and how to apply them. From a CSE standpoint, future work could identify the unique needs of mathematical modelers across different levels of real-world variability and system adaptation needs. This includes tailoring knowledge elicitation to the specific categories of variability and the corresponding model adaptations outlined in Sharkey et al. (2024).

Throughout the case study, each of the model iterations offers a way to address the decision of where to locate hand sanitizer stations. Besides academic buildings, the model can also accommodate residence halls (i.e., on-campus student living) and dining halls as other possible locations. These types of buildings account for significant foot traffic on campus as well. While the final model considers locations near exterior doors only, another area of interest that could be addressed could be building interiors. The limiting factor is data availability, but if other sources become available, e.g., WiFi access points, another iteration could be conducted to reflect this change. The models could also be used to locate other demand-based items at the university, including garbage cans, marketing materials for public health strategies, etc. Additionally, the presented location models could be applied at other universities. Many schools similarly use card access in the same manner to gain access into buildings. In these cases, data for the specific school could be incorporated into the final model version to distribute hand sanitizers or other demand-based resources.

## 6. Conclusions

This paper seeks to improve the quality of optimization model recommendations through a rigorous integration of stakeholders into the modeling process. CSE methods provide a framework to ground modelers' efforts to understand novel and evolving situations. The developed models reflect the decision-maker strategies, priorities, and constraints as well as input from other key stakeholders. The integrated OR-CSE approach is broadly applicable in other settings. The process may be particularly useful when both stakeholders' and modelers' understanding of situations change over time, e.g., changing information related to COVID-19. By leveraging the quantitative and qualitative context, we suggest that this process may enhance the effectiveness of recommendations produced by optimization models.

## References

Ahani, N., & Trapp, A. C. (2021). *Human-Centric Decision Support Tools: Insights from Real-World Design and Implementation*. https://arxiv.org/ftp/arxiv/papers/2111/2111.05796.pdf

Andrews, J. M., Farias, V. F., Khojandi, A. I., & Yan, C. M. (2019). Primal–Dual Algorithms for Order Fulfillment at Urban Outfitters, Inc. *INFORMS Journal on Applied Analytics*, *49*(5), 355–370. https://doi.org/10.1287/inte.2019.1013

Arnstein, S. R. (1969). A Ladder Of Citizen Participation. *Journal of the American Planning Association*, *35*(4), 216–224. https://doi.org/10.1080/01944366908977225

Barnhart, C., Belobaba, P., & Odoni, A. R. (2003). Applications of operations research in the air transport industry. *Transportation Science*, *37*(4), 368–391. https://doi.org/10.1287/trsc.37.4.368.23276

Barnhart, C., Bertsimas, D., Delarue, A., & Yan, J. (2022). Course Scheduling Under Sudden Scarcity: Applications to Pandemic Planning. *Manufacturing and Service Operations Management*, *24*(2), 727–745. https://doi.org/10.1287/msom.2021.0996

Bennett, K. B., Edman, C., Cravens, D., & Jackson, N. (2023). Decision Support for Flexible Manufacturing Systems: Application of the Cognitive Systems Engineering and Ecological Interface Design Approach. *Journal of Cognitive Engineering and Decision Making*, *17*(2), 99–119. https://doi.org/10.1177/15553434221118976

Birmingham, W. C., Wadsworth, L. L., Lassetter, J. H., Graff, T. C., Lauren, E., & Hung, M. (2023). COVID-19 lockdown: Impact on college students' lives. *Journal of American College Health*, *71*(3), 879–893. https://doi.org/10.1080/07448481.2021.1909041

Bisantz, A. M., & Burns, C. M. (Eds.). (2009). *Applications of Cognitive Work Analysis*. CRC Press.

Bojke, L., Claxton, K., Sculpher, M., & Palmer, S. (2009). Characterizing structural uncertainty in decision analytic models: A review and application of methods. *Value in Health*, *12*(5), 739–749. https://doi.org/10.1111/j.1524-4733.2008.00502.x

Braun, V., & Clarke, V. (2012). Thematic analysis. In *APA handbook of research methods in psychology, Vol 2: Research designs: Quantitative, qualitative, neuropsychological, and biological.* (pp. 57–71). American Psychological Association. https://doi.org/10.1037/13620-004

Burns, C. M., & Hajdukiewicz, J. (2017). *Ecological Interface Design*. CRC Press. https://doi.org/10.1201/9781315272665

Centers for Disease Control and Prevention (CDC). (2020). *Cleaning, Disinfection, and Hand Hygiene in Schools – a Toolkit for School Administrators*. Coronavirus Disease 2019 (COVID-19). https://stacks.cdc.gov/view/cdc/97530

Chen, G., Fei, X., Jia, H., Yu, X., & Shen, S. (2022). The University of Michigan Implements a Hub-and-Spoke Design to Accommodate Social Distancing in the Campus Bus System Under COVID-19 Restrictions. *INFORMS Journal on Applied Analytics*, *52*(6), 539–552. https://doi.org/10.1287/inte.2022.1131

Chung, Q. B., Willemain, T. R., & O'Keefe, R. M. (2000). Influence of model management systems on decision making: Empirical evidence and implications. *Journal of the Operational Research Society*, *51*(8), 936–948. https://doi.org/10.1057/palgrave.jors.2600993

Currie, C. S. M., Fowler, J. W., Kotiadis, K., Monks, T., Onggo, B. S., Robertson, D. A., & Tako, A. A. (2020). How simulation modelling can help reduce the impact of COVID-19. *Journal of Simulation*, *14*(2), 83–97. https://doi.org/10.1080/17477778.2020.1751570

Daskin, M. S., & Tucker, E. L. (2018). The trade-off between the median and range of assigned demand in facility location models. *International Journal of Production Research*, *56*(1–2), 97–119. https://doi.org/10.1080/00207543.2017.1401751

de Gooyert, V., Rouwette, E., van Kranenburg, H., & Freeman, E. (2017). Reviewing the role of stakeholders in Operational Research: A stakeholder theory perspective. *European Journal of Operational Research*, *262*(2), 402–410. https://doi.org/10.1016/j.ejor.2017.03.079

Doyle, J. K., & Ford, D. N. (1998). Mental models concepts for system dynamics research. *System Dynamics Review*, *14*(1), 3–29. https://doi.org/10.1002/(SICI)1099-1727(199821)14:1<3::AID-SDR140>3.0.CO;2-K

Durán, G., Guajardo, M., Gutiérrez, F., Marenco, J., Sauré, D., & Zamorano, G. (2021). Scheduling the Main Professional Football League of Argentina. *INFORMS Journal on Applied Analytics*, *51*(5), 361–372. https://doi.org/10.1287/inte.2021.1088

Duran, S., Gutierrez, M. A., & Keskinocak, P. (2011). Pre-Positioning of Emergency Items for CARE International. *Interfaces*, *41*(3). https://doi.org/10.1287/inte.ll00.0526

Eddy, D. M., Hollingworth, W., Caro, J. J., Tsevat, J., McDonald, K. M., & Wong, J. B. (2012). Model transparency and validation: A report of the ISPOR-SMDM modeling good research practices task force-7. *Medical Decision Making*, *32*(5), 733–743. https://doi.org/10.1177/0272989X12454579

Eisenberg, D., Seager, T., & Alderson, D. L. (2019). Rethinking Resilience Analytics. *Risk Analysis*, *39*(9), 1870–1885. https://doi.org/10.1111/risa.13328

Elix, B., & Naikar, N. (2021). Designing for Adaptation in Workers' Individual Behaviors and Collective Structures With Cognitive Work Analysis: Case Study of the Diagram of Work Organization Possibilities. *Human Factors*, *63*(2), 274–295.

https://doi.org/10.1177/0018720819893510

Elm, W. C., Gualtieri, J. W., Mckenna, B. P., Tittle, J. S., Peffer, J. E., Szymczak, S. S., & Grossman, J. B. (2008). Integrating Cognitive Systems Engineering Throughout the Systems Engineering Process. *Journal of Cognitive Engineering and Decision Making*, *2*(3), 249–273. https://doi.org/10.1518/155534308X377108

Facilities Support Services. (2020). *Freeman Hall*.

Foster, S. A., Hegde, S., O'Brien, T. C., & Tucker, E. L. (2023). Organizational Adaptive Capacity During a Large-Scale Surprise Event: A Case Study at an Academic Institution During the COVID-19 Pandemic. *IISE Transactions on Occupational Ergonomics and Human Factors*, *0*(0), 1–20. https://doi.org/10.1080/24725838.2023.2221045

Gass, S. I. (1983). Decision-Aiding Models: Validation, Assessment, and Related Issues for Policy Analysis. *Operations Research*, *31*(4), 603–631. https://doi.org/10.1287/opre.31.4.603

Gentry, S., Chow, E., Massie, A., & Segev, D. (2015). Gerrymandering for Justice: Redistricting U.S. Liver Allocation. *INFORMS Journal on Applied Analytics*, *45*(5), 462–480. https://doi.org/10.1287/inte.2018.0980

Geoffrion, A. M. (1976). The Purpose of Mathematical Programming is Insight, Not Numbers. *Interfaces*, *7*(1), 81–92. https://doi.org/10.1287/inte.7.1.81

Goncalves Filho, A. P., Jun, G. T., & Waterson, P. (2019). Four studies, two methods, one accident – An examination of the reliability and validity of Accimap and STAMP for accident analysis. *Safety Science*, *113*(October 2018), 310–317. https://doi.org/10.1016/j.ssci.2018.12.002

Gore, A. B., Kurz, M. E., Saltzman, M. J., Splitter, B., Bridges, W. C., & Calkin, N. J. (2022).

Clemson University's Rotational Attendance Plan During COVID-19. *INFORMS Journal on Applied Analytics*, *52*(6), 553–567. https://doi.org/10.1287/inte.2022.1139

Hegde, S., Hettinger, A. Z., Fairbanks, R. J., Wreathall, J., Krevat, S. A., & Bisantz, A. M. (2020). Knowledge Elicitation to Understand Resilience: A Method and Findings From a Health Care Case Study. *Journal of Cognitive Engineering and Decision Making*, *14*(1), 75–95. https://doi.org/10.1177/1555343419877719

Hettinger, A. Z., Roth, E. M., & Bisantz, A. M. (2017). Cognitive engineering and health informatics: Applications and intersections. *Journal of Biomedical Informatics*, *67*, 21–33. https://doi.org/10.1016/j.jbi.2017.01.010

Hollnagel, E., & Woods, D. D. (2005). *Joint Cognitive Systems: Foundations of Cognitive Systems Engineering* (1st ed.). CRC Press. https://doi.org/10.1201/9781420038194

Hollnagel, E., Woods, D. D., & Leveson, N. (2006). *Resilience Engineering: Concepts and Precepts*. Ashgate Publishing Ltd.

Johnson, M. P., Midgley, G., & Chichirau, G. (2018). Emerging trends and new frontiers in community operational research. *European Journal of Operational Research*, *268*(3), 1178–1191. https://doi.org/10.1016/j.ejor.2017.11.032

Jones, L., & Wells, K. (2007). Strategies for Academic and Clinician Engagement in Community-Participatory Partnered Research. *JAMA*, *297*(4), 407. https://doi.org/10.1001/jama.297.4.407

Klein, G. a, Calderwood, R., & Macgregor, D. (1989). Critical Decision Method for Eliciting Knowledge. *IEEE Transactions On Systems Man And Cybernetics*, *19*(3), 462–472. http://ieeexplore.ieee.org/xpls/abs_all.jsp?arnumber=31053

Kovacs, G., & Moshtari, M. (2019). A roadmap for higher research quality in humanitarian operations: A methodological perspective. *European Journal of Operational Research*, *276*(2), 395–408. https://doi.org/10.1016/j.ejor.2018.07.052

Kunz, N., Van Wassenhove, L. N., Besiou, M., Hambye, C., & Kovács, G. (2017). Relevance of humanitarian logistics research: best practices and way forward. *International Journal of Operations and Production Management*, *37*(11), 1585–1599. https://doi.org/10.1108/IJOPM-04-2016-0202

Mayring, P. (2004). Qualitative Content Analysis. *A Companion to Qualitative Research*, *1*(2), 159–176.

McGeorge, N., Hegde, S., Berg, R. L., Guarrera-Schick, T. K., LaVergne, D. T., Casucci, S. N., Hettinger, A. Z., Clark, L. N., Lin, L., Fairbanks, R. J., Benda, N. C., Sun, L., Wears, R. L., Perry, S., & Bisantz, A. (2015). Assessment of innovative emergency department information displays in a clinical simulation center. *Journal of Cognitive Engineering and Decision Making*, *9*(4), 329–346. https://doi.org/10.1177/1555343415613723

Merrick, J. H., & Weyant, J. P. (2019). On choosing the resolution of normative models. *European Journal of Operational Research*, *279*(2), 511–523. https://doi.org/10.1016/j.ejor.2019.06.017

Miles, M. B., & Huberman, A. M. (1994). *Qualitative data analysis: An expanded sourcebook* (2nd ed.). SAGE Publications, Inc.

Militello, L. G., Dominguez, C. O., Lintern, G., & Klein, G. (2010). The role of cognitive systems engineering in the systems engineering design process. *Systems Engineering*, *13*(3), 261–273. https://doi.org/10.1002/sys.20147

Mingers, J. (2011). Soft OR comes of age—but not everywhere! *Omega*, *39*(6), 729–741. https://doi.org/10.1016/j.omega.2011.01.005

Mingers, J., & Rosenhead, J. (2004). Problem structuring methods in action. *European Journal of Operational Research*, *152*(3), 530–554. https://doi.org/10.1016/S0377-2217(03)00056-0

Morgan, K., Collier, Z. A., Gilmore, E., & Schmitt, K. (2022). Decision-first modeling should guide decision making for emerging risks. *Risk Analysis*, *42*(12), 2613–2619. https://doi.org/10.1111/risa.13888

Morris, W. T. (1967). On the art of modeling. *Management Science*, *13*(12), B707–B717.

Murphy, F. H. (2005). ASP, The Art and Science of Practice: Elements of a Theory of the Practice of Operations Research: A Framework. *Interfaces*, *35*(2), 154–163.

Navabi-Shirazi, M., El Tonbari, M., Boland, N., Nazzal, D., & Steimle, L. N. (2022). Multicriteria Course Mode Selection and Classroom Assignment Under Sudden Space Scarcity. *Manufacturing and Service Operations Management*, *24*(6), 3252–3268. https://doi.org/10.1287/msom.2022.1131

O'Brien, T., Foster, S., Tucker, E. L., & Hegde, S. (2021). COVID Response: A Blended Approach to Studying Sanitizer Station Deployment at a Large Public University. *2021 Resilience Week (RWS)*.

Pidd, M. (2010). Why modelling and model use matter. *Journal of the Operational Research Society*, *61*(1), 14–24. https://doi.org/10.1057/jors.2009.141

Pidd, Michael. (1999). Just modeling through: A rough guide to modeling. *Interfaces*, *29*(2), 118–132. https://doi.org/10.1287/inte.29.2.118

Pollack, M., & Steimle, L. N. (2022). The Implications of State Aggregation on the Utility of Estimated Markov Decision Process Models. *Preprint*.

Rönnqvist, M., Svenson, G., Flisberg, P., & Jönsson, L. E. (2017). Calibrated route finder: Improving the safety, environmental consciousness, and cost effectiveness of truck routing

in Sweden. *Interfaces*, *47*(5), 372–395. https://doi.org/10.1287/inte.2017.0906

Sharkey, T. C., Foster, S., Hegde, S., Kurz, M. E., & Tucker, E. L. (2024). A categorization of observed uses of operational research models for fundamental surprise events: Observations from university operations during COVID-19. *Journal of the Operational Research Society*, 1–13. https://doi.org/10.1080/01605682.2024.2346117

Sharkey, T. C., Nurre Pinkley, S. G., Eisenberg, D. A., & Alderson, D. L. (2021). In search of network resilience: An optimization-based view. *Networks*, *77*(2), 225–254. https://doi.org/10.1002/net.21996

Shehadeh, K. S., & Tucker, E. L. (2022). Stochastic optimization models for location and inventory prepositioning of disaster relief supplies. *Transportation Research Part C: Emerging Technologies*, *144*(August), 103871. https://doi.org/10.1016/j.trc.2022.103871

Shen, S. (2020). From Data to Actions, From Observations to Solutions: A Summary of Operations Research and Industrial Engineering Tools for Fighting COVID-19. *Report*. https://public.websites.umich.edu/~siqian/docs/presentation/or-ie-fighting-covid19_v1.pdf

Sodhi, M. S., & Tang, C. S. (2014). Guiding the next generation of doctoral students in operations management. *International Journal of Production Economics*, *150*, 28–36. https://doi.org/10.1016/j.ijpe.2013.11.016

Spetzler, C. S., & Stael von Holstein, C. A. S. (1975). Probability Encoding in Decision Analysis. *Management Science*, *22*(3), 340–358. https://doi.org/10.1287/mnsc.22.3.340

Steimle, L. N., Sun, Y., Johnson, L., Besedeš, T., Mokhtarian, P., & Nazzal, D. (2022). Students' preferences for returning to colleges and universities during the COVID-19 pandemic: A discrete choice experiment. *Socio-Economic Planning Sciences*, *82*(February). https://doi.org/10.1016/j.seps.2022.101266

Strong, M., Oakley, J. E., & Chilcott, J. (2012). Managing structural uncertainty in health economic decision models: A discrepancy approach. *Journal of the Royal Statistical Society. Series C: Applied Statistics*, *61*(1), 25–45. https://doi.org/10.1111/j.1467-9876.2011.01014.x

Willemain, T. R. (1994). Insights on Modeling from a Dozen Experts. *Operations Research*, *42*(2), 213–222. https://doi.org/10.1287/opre.42.2.213

Willemain, T. R. (1995). Model Formulation: What Experts Think about and When. *Operations Research*, *43*(6), 916–932.

Willemain, T., Wallace, W., Felichmann, K., Waisel, L., & Ganaway, S. (2003). Bad numbers: coping with flawed decision support. *Journal of the Operational Research Society*, *54*(9), 949–957.

Woods, D. D., & Roth, E. M. (1988). Cognitive Systems Engineering. In *Handbook of Human-Computer Interaction* (pp. 3–43). North Holland. https://doi.org/10.1016/B978-0-444-70536-5.50006-3

Appendix to
"Developing Optimization Models with Cognitive Systems Engineering"

The semi-structured interview protocol is provided in Section A-1. Additional figures are presented in Section A-2, and additional model results are shown in Section A-3.

## Section A-1. Semi-Structured Interview Protocol

Brief introduction about the study and purpose of the interview. *"The purpose of this study is to understand how Clemson University adapted to the COVID-19 pandemic. We are particularly interested in understanding how decisions were made relating to hand sanitation strategies on campus and how these strategies have been implemented. During this interview, we will ask you to describe your role in Clemson's response and hand sanitation policies, and seek your perspectives on what worked and what didn't. We seek your permission to continue recording this interview for the purpose of transcription and analysis. Do you have any questions or concerns before we begin?"*

<u>Note to interviewer:</u> The questions listed below represent important themes or categories of information to be obtained from the participant. The order of questions can be adjusted according to the flow of the interview and topics emerging from the participant's responses. 'Probes' listed alongside some of the questions indicate items or themes that the interviewer should look out for in the responses and ensure are covered during the interview. They may be used as follow up questions to guide the response of the participant or to have them elaborate on a response.

### Questions

Briefly describe your role at Clemson University and how it relates to the University's hand sanitation policies.

What were the main challenges faced at the start of the pandemic (around March 2020) in terms of decisions around facilities / hand sanitation?

Was there a contingency plan available for such a situation?

What types of information were needed or sought regarding hand sanitation? Were these available?

> <u>Probes:</u> Categories and sources of information

How was this information shared?

> <u>Probes</u>: Media of communication – emails, dashboards, databases, messaging systems

How did you know what information to look at and when?

> <u>Probes</u>: Government and external reports and guidelines (e.g. CDC, WHO, Federal, SC State)

Did you consider students' and their families' perspectives during this process? How?

What were some of the considerations that emerged?

> <u>Probes</u>: Safety; flexibility; accommodation (esp. for those coming from outside of Clemson)

What were your main priorities in terms of hand sanitation?

Probes: Infection control; graduation times; student stress; resources for remote classes

Did you use any data analytics or modeling in the station placement? Please describe this process.

How did you decide what might be an 'optimum' number of stations for individual buildings and the entire campus? Was there a prioritization of building coverage over campus coverage?

Probes: no. of classrooms or spaces; classroom size; no. of students who wanted to be on campus vs. remote; other

How did variability and uncertainty influence these methods?

Probes: changes in advisories/guidelines; infection rates; student attendance; other

Was there any simulation or testing involved?

Once classes resumed (hybrid / full in-person), was there a way to monitor how things were going?

Probes: attendance; infection rates; testing protocols; compliance; other

Did plans/strategies have to be adjusted based on actual dynamics of classes and student-presence on campus? What were some of these changes?

Looking back, what types of information/data have been most important or useful in this process? What additional data or information would have been useful?

How was information and decisions communicated across institutional layers – President and Provost's offices; college-level leadership; departmental leadership; other administrative offices?

Probes: Emails, townhalls, dashboards, meetings/taskforce

How were decisions and implementation actions coordinated across various bodies? Who were the facilitators?

What directions are being looked at for the future? Will strategies drastically change?

What are the current benchmark goals for right now and future semesters?

## Section A-2. Additional Figures

Figure A-1 is a picture of a hand sanitizer station. Figure A-2 presents a heat map of door accesses by building in February 2021. It was produced using Mapline (Mapline, 2022) and presented with permission.

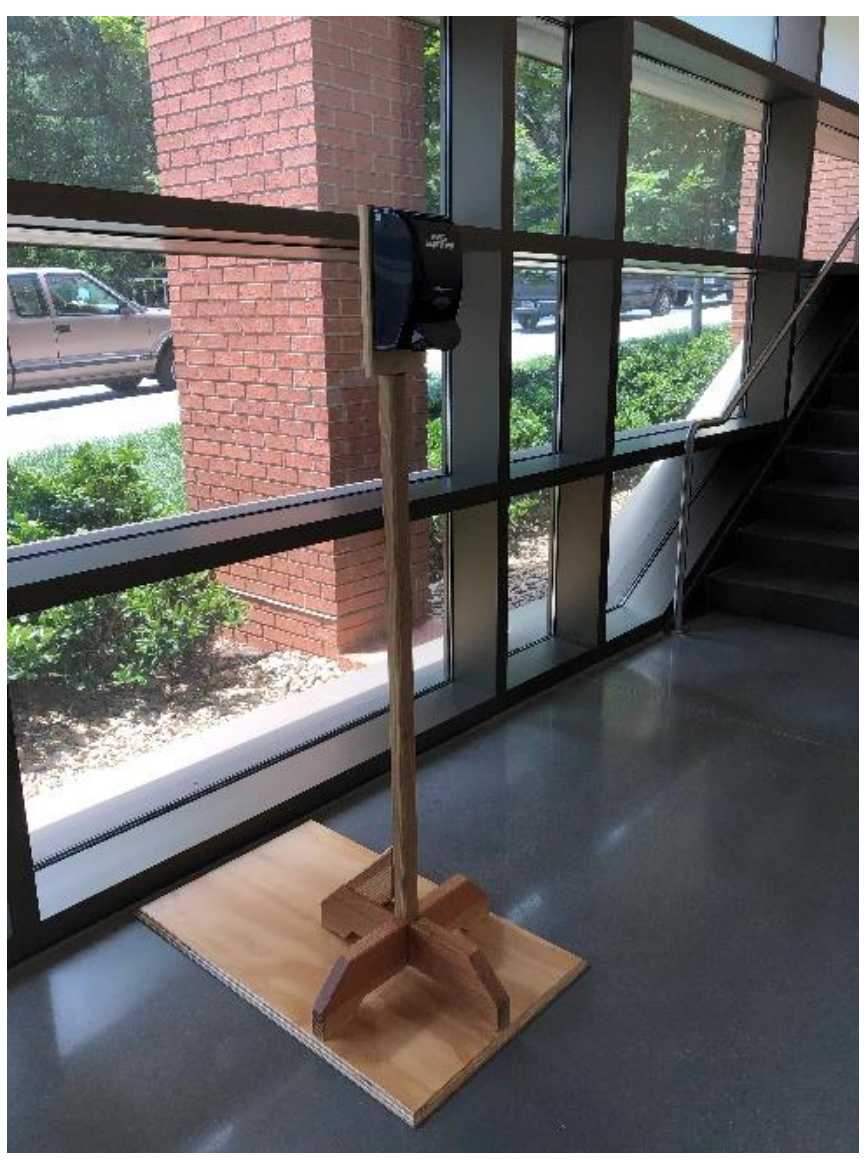

**Figure A-6.** Hand sanitizer

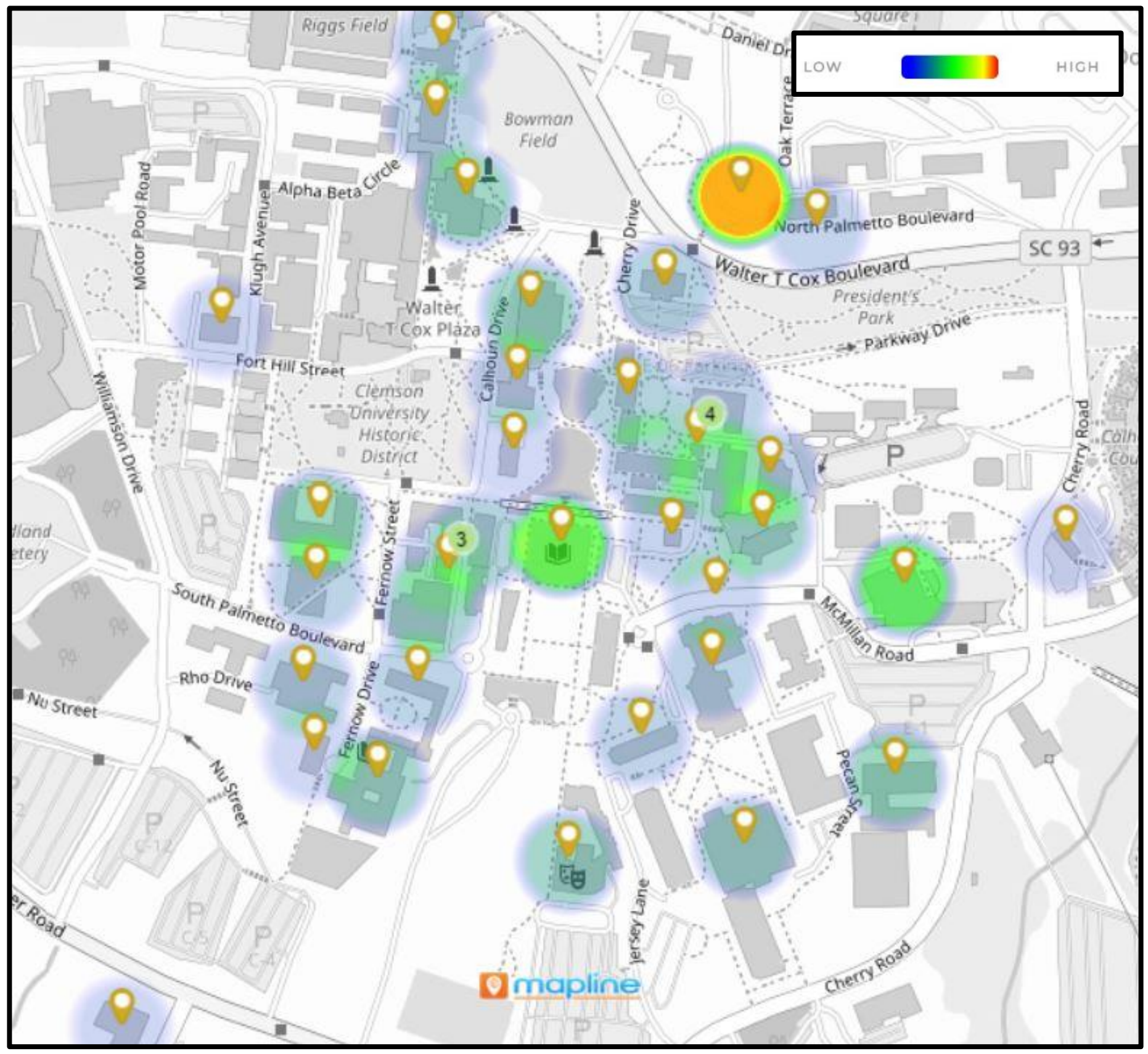

**Figure A-2.** Heat map of door access control data by building in Feb.

### Section A-3. Additional Modeling Results

**Table A-1:** Number of dispensers allocated to each building

| Building Name | Initial Deployment | V1: Facilities Heuristic | V3: Target Coverage | V5: Max Coverage by Door |
|---|---|---|---|---|
| Academic Success Center | 2 | 2 | 1 | 2 |
| Administrative Services Building | 1 | 1 | 1 | 7 |
| Barre Hall | 4 | 2 | 1 | 4 |
| BioSystems Research Center | 4 | 5 | 4 | 5 |
| Brackett Hall | 2 | 4 | 4 | 4 |
| Brooks Center | 2 | 4 | 3 | 3 |
| Campbell Museum | 1 | 1 | 1 | 1 |
| College of Business Building | 20 | 15 | 24 | 9 |
| Cook Laboratory | 1 | 1 | 1 | 1 |
| Dillard Building | 2 | 1 | 1 | 1 |
| Earle Hall | 1 | 2 | 1 | 3 |
| Edwards Hall | 2 | 4 | 3 | 4 |
| Fluor Daniel | 1 | 2 | 1 | 3 |
| Freeman Hall | 3 | 4 | 2 | 4 |
| Godfrey Hall | 1 | 2 | 1 | 1 |
| Hardin Hall | 3 | 2 | 1 | 1 |
| Harris Smith | 1 | 1 | 1 | 1 |
| Holtzendorff Hall | 1 | 3 | 1 | 1 |
| Hunter Hall | 1 | 3 | 1 | 3 |
| Jordan Hall | 2 | 2 | 1 | 1 |
| Kinard Hall | 4 | 2 | 1 | 2 |
| Lee Hall (3 Buildings) | 6 | 4 | 3 | 2 |
| Long Hall | 1 | 2 | 1 | 2 |
| Lowry Hall | 1 | 2 | 1 | 2 |
| Martin Hall (3 Buildings) | 6 | 3 | 1 | 3 |
| McAdams Hall | 2 | 3 | 1 | 2 |
| Old Main | 6 | 5 | 4 | 4 |
| Olin Hall | 1 | 1 | 1 | 1 |
| P&A Building | 4 | 4 | 3 | 4 |
| Rhodes Hall/Annex | 2 | 3 | 1 | 2 |
| Cooper Library | 4 | 9 | 13 | 4 |
| Sikes Hall | 3 | 3 | 1 | 3 |
| Sirrine Hall | 4 | 5 | 4 | 5 |
| Strode Tower | 1 | 2 | 1 | 2 |

| Vickery Hall | 1 | 2 | 1 | 1 |
|---|---|---|---|---|
| Watt Innovation | 1 | 8 | 11 | 4 |

Table A-2 presents additional results for Model V2. In this model, the definition of a covered classroom depends on how long people are willing to walk to use a station. We define a coverage time threshold that indicates the time between a station and the furthest classroom. With 3 or 4 stations, the coverage time is 8 seconds. With 1 or 2 stations, the time increases substantially to 30 or 24 seconds, respectively.

**Table A-2**. Optimal locations and maximum threshold for coverage time

| **Dispensers allocation (*p*)** | **Optimal Dispenser Locations** | **Maximum Coverage Time** |
|---|---|---|
| $p = 4$ | 1, 5, 9, 10 | $T = 8$ sec |
| $p = 3$ | 5, 9, 10 | $T = 8$ sec |
| $p = 2$ | 5, 10 | $T = 24$ sec |
| $p = 1$ | 1 | $T = 30$ sec |

**Reference**

Mapline Inc. (2022). *Mapline*. https://mapline.com/